\documentclass[10pt,a4paper]{article}
\usepackage[utf8]{inputenc}
\usepackage{amsmath}
\usepackage{amsfonts}
\usepackage{amssymb}
\usepackage{amsthm}
\newtheorem{theorem}{Theorem}

\newtheorem*{theorem-a}{Theorem A}
\newtheorem*{theorem-b}{Theorem B}
\newtheorem{definition}{Definition}

\usepackage{hyperref}
\newcommand{\footremember}[2]{%
    \footnote{#2}
    \newcounter{#1}
    \setcounter{#1}{\value{footnote}}%
}
 
\title{Nonfreeness of some algebras of hermitian modular forms.}
\author{%
  Stuken Ekaterina\footremember{alley}{Chebyshev Laboratory, St. Petersburg State University, 14th Line V.O., 29B, Saint Petersburg 199178
Russia.
}%
  }

\begin{document}
\maketitle
\begin{abstract}
We study the algebras of hermitian automorphic forms for the lattice $L_n=diag(1,1,\ldots,1,-1)$ and for the field $K=\mathbb{Q}(\sqrt{-d})$ such that $p=2$ is unramified and the ring of integers $\mathcal{O}_K$ is a p.i.d. We prove that for $d>7$ these algebras can't be free. When $d=7$ and $d=3$ we give an estimate for the dimension of the symmetric spaces for which these algebras might be free. We also compare our results with the known results for $d=3$.
\end{abstract}
\section*{Introduction}

Let $K=\mathbb{Q}(\sqrt{-d})$ be an imaginary quadratic field, $d$ is odd and square free. We restrict ourselves to the case when the ring of integers $\mathcal{O}_K$ is a principal ideal domain and $p=2$ is unramified. 
Let $L$ be an integral hermitian lattice of signature $(n,1),~n>1$. Let $V$ be the vector space $V=L\otimes_{\mathcal{O}_K} \mathbb{C}$. Denote by $\Gamma$ the group $U(L)$. We consider the  Hermitian symmetric domain $\mathcal{D} = \{[z]\in P(V) | (z, z) < 0\}$ and a principal $\mathbb{C}^*$-bundle $\mathcal{A}=\{z\in V | [z]\in \mathcal{D}\}$.

% Let $V$ be the vector space $V=L_n\otimes_{\mathcal{O}_K} \mathbb{C}$. We consider the  Hermitian symmetric domain $\mathcal{D} = \{[z]\in P(V) | (z, z) < 0\}$ and a principal $\mathbb{C}^*$-bundle $\mathcal{A}=\{z\in V | [z]\in \mathcal{D}\}$.

\begin{definition} Automorphic form of weight $k\in \mathbb{Z}, ~k>0$ for $\Gamma$ with character $\chi:\Gamma\rightarrow \mathbb{C}^*$ is a holomorphic function $f$ in $\mathcal{A}$ such that:\\
1) $f(tz)=t^{-k}f(z)$, $t\in\mathbb{C}^*$\\
2) $f(g(z))=\chi(g)f(z)$, $g\in \Gamma$. 
\end{definition}

%\begin{definition} The Hirzebruch-Mumford volume of the orthogonal group $O^+(L)$ of the lattice $L$ is a normalised covolume $\mathrm{Vol}_{HM}(D/O^+(L))=\frac{\mathrm{Vol}(D/O^+(L))}{\mathrm{Vol}(D^c)}$, where $D^c$ is a compact Hermitian symmetric space dual to the space $D$.
%\end{definition}

%The question of whether the algebra $A(\Gamma)$ is free is one of the main questions in transcendental invariant theory (\cite{vin2}). The purpose of this paper is to find all $d$, for which the algebra $A(\Gamma)$ can be free. The following theorem gives the answer.
Here we prove the following theorem:
\begin{theorem}\label{main}
Let $K=\mathbb{Q}(\sqrt{-d})$ be a quadratic field such that it's ring of integers is a p.i.d. and $p=2$ is unramified. Consider the lattice $L_n=diag(1,1,\ldots,1,-1)$ of signature $(n,1)$. Let $A_n$ be the algebra of automorphic forms for the group $U(L_n)$. Then:
\begin{itemize}
\item if $d>7$ than the algebra $A_n$ can't be free;
\item if $d=7$ than the algebra $A_n$ can be free only for $n\leq 4$;
\item if $d=3$ than the the algebra $A_n$ can be free only for $n\leq 7$.
\end{itemize}
\end{theorem}

To prove this theorem we consider the invariant $K(\Gamma)$ for the group $\Gamma$:
\begin{equation}
K(\Gamma):=\dfrac{\sum\limits_{[\pi]} \frac{n_\pi-1}{n_\pi}\mathrm{Vol}_{HM}(B^{n-1}/\Gamma_\pi)}{\mathrm{Vol}_{HM}(B^n/\Gamma)},
\end{equation}
where $\Gamma_\pi$ is the stabiliser of the mirror $\pi$ in the group $\Gamma$, the sum is taken over all $\Gamma$-conjugacy classes of the mirrors of reflections $\pi$ in the group $\Gamma$, $n_\pi$ is the order of the reflection in $\pi$.

We will need one theorem of J.\,H.~Bruinier (\cite{bru}). In the case of interest to us and in a form convenient for us it states the following: 
\begin{theorem-a}
Let $F$ be a meromorphic automorphic form for the group $\Gamma$ of weight $K$ such that it has zeros on all mirrors of reflections in $\Gamma$ and only on them. Denote by $n_\pi$ the order of the reflection in $\pi$. If each zero has the order $n_\pi-1$, then $K(\Gamma)=K$.
\end{theorem-a}

We also use another theorem of H.~Aoki, T.~Ibukiyama (\cite{aok}):
\begin{theorem-b} If the algebra of $\Gamma$-automorphic forms is free with the generators of weights $k_1,~k_2,\ldots, ~k_{n+1}$, then there exists unique up to the proportionality $\Gamma$-automorphic form $F$ of weight $n+1+k_1+k_2+\ldots+k_{n+1}$ whose divisor is $div(F)=\sum (n_\pi-1)\pi$ (the sum is taken over all mirrors of reflections in the group $\Gamma$).
\end{theorem-b}
We apply this theorem to the lattice $L_n$ and the group $\Gamma=U(L_n)$. If the algebra $A_n$ is free then it has $n+1$ generators of weight at least $1$. Hence the weight $K$ of the form $F$ from theorem $B$ for $\Gamma$ is not less then $2n+2$. Applying to the form $F$ theorem $A$ we get that $K(\Gamma)\geq 2n+2$. 

The idea of the proof of theorem \ref{main} is to estimate in formula (1) the value of $K(\Gamma)$. We show that the denominator of $K(\Gamma)$ grows much faster than the numerator, thus $K(\Gamma)$ is less then $2n+2$ when $n$ is great enough. %numerator from the above and the denominator from the below. We show that the denominator grows faster as a function of $n$, thus $K(\Gamma)$ is less then $2n+2$ when $n$ is great enough. 

%We consider the invariant $K(\Gamma)$ for the group $\Gamma$, which is a sum over classes of $\Gamma$-conjugated roots of the lattice $L_n$ of covolumes of stabilisers of subballs $B^{n-1}$, divided by the covolume of the lattice $L_n$ (\cite{bru}, \cite{stu}). 
%In this article we estimate it from the above and use this estimate to prove that the algebra of $\Gamma$-automorphic forms can't be free for all but a finite number of values of $n$.

  H.Wang and B.Williams in \cite{wan} proved that for $d=3$ the algebra $A_n$ is free when $n\leq 4$ and computed the weights of generators. One can easily check that these weights correspond to our calculations of the value $K(\Gamma)$, i.e. $K(\Gamma)=K=n+1+k_1+k_2+\ldots+k_{n+1}$.

Let $L_n=\mathrm{diag}(1,1,\ldots,1,-1)$ and $M_n=\mathrm{diag}(1,1,\ldots,1,-2)$ be the lattices of signature $(n,1)$. Consider the groups $SU(L_n,\mathcal{O}_K)$ and $SU(M_n,\mathcal{O}_K)$ for the ring of integers $\mathcal{O}_K$ of the field $K$. These are discrete arithmetic groups naturally acting on a complex ball $B^n$. Their quotient spaces have finite volumes. We will write $\mathrm{Vol}(L_n)$ and $\mathrm{Vol}(M_n)$ respectively for the Hirzebruch-Mumford covolume of the spaces $\mathbb{B}^n/SU(L_n,\mathcal{O}_K)$ and $\mathbb{B}^n/SU(M_n,\mathcal{O}_K)$.

To prove theorem \ref{main} we need to know the values of covolumes $\mathrm{Vol}(L_n)$ and $\mathrm{Vol}(M_n)$. The first one was calculated in \cite{zel}, the second in \cite{stu}. Denote by $L(k)$ the $L$-function of the quadratic field $K$ with character $\chi_D(p)=\left(\frac{D}{p}\right)$ (Kronecker symbol).
\newpage
\begin{theorem}\label{L}
The Hirzebruch-Mumford covolume of $\mathbb{B}^n/SU(L_n,\mathcal{O}_K)$ is
\begin{table}[h]
\centering
	\begin{tabular}{|c|c|c|}
		\hline
		{\bf n} & {\bf D} & {\bf $\mathrm{Vol}(L_n)$} \\
	 	\hline
	 	{\bf even} &  {\bf -4d } & {\bf $D^{\frac{n^2+3n}{4}}\cdot \prod\limits_{j=1}^{n}\frac{j!}{(2\pi)^{j+1}}\cdot \zeta(2)\cdot L(3)\cdot \zeta(4)\cdot L(5)\cdot \ldots \cdot L(n+1)$}   \\
	 	\hline
	 	{\bf even} &  {\bf -d } & {\bf $D^{\frac{n^2+3n}{4}}\cdot \prod\limits_{j=1}^{n}\frac{j!}{(2\pi)^{j+1}}\cdot \zeta(2)\cdot L(3)\cdot \zeta(4)\cdot L(5)\cdot \ldots \cdot L(n+1)$}   \\
	 	\hline
	 	{\bf odd} &  {\bf -4d} & {\bf $D^{\frac{n^2+3n}{4}}\cdot(1-2^{-(n+1)})\cdot \prod\limits_{p|d} (1+\left( \frac{(-1)^{(n+3)/2}}{p}\right)\cdot p^{-\frac{n+1}{2}}) \prod\limits_{j=1}^{n}\frac{j!}{(2\pi)^{j+1}}\cdot \zeta(2)\cdot  \ldots \cdot \zeta(n+1)$}   \\
	 	\hline
	 	{\bf odd} &  {\bf -d} & {\bf $D^{\frac{n^2+3n}{4}}\cdot \prod\limits_{p|d} (1+\left( \frac{(-1)^{(n+3)/2}}{p}\right)p^{-\frac{n+1}{2}}) \prod\limits_{j=1}^{n}\frac{j!}{(2\pi)^{j+1}}\cdot \zeta(2)\cdot L(3)\cdot  \ldots \cdot \zeta(n+1)$}      \\
	 	\hline
	\end{tabular}
\end{table}
\end{theorem}
%\vspace{-0.5cm}

\begin{theorem}\label{M}
The Hirzebruch-Mumford covolume of $\mathbb{B}^n/SU(M_n,\mathcal{O}_K)$ is \begin{table}[h]
\centering
	\begin{tabular}{|c|c|c|}
		\hline
		{\bf n} & {\bf D} & {\bf $\mathrm{Vol}(M_n)$} \\
	 	\hline
	 	{\bf even} &  {\bf -4d } & {\bf $D^{\frac{n^2+3n}{4}}\cdot \prod\limits_{j=1}^{n}\frac{j!}{(2\pi)^{j+1}}\cdot \zeta(2)\cdot L(3)\cdot \zeta(4)\cdot L(5)\cdot \ldots \cdot L(n+1)\cdot (2^n-1)$}   \\
	 	\hline
	 	{\bf even} &  {\bf -d} & {\bf $D^{\frac{n^2+3n}{4}}\cdot \prod\limits_{j=1}^{n}\frac{j!}{(2\pi)^{j+1}}\cdot \zeta(2)\cdot L(3)\cdot \zeta(4)\cdot L(5)\cdot \ldots \cdot L(n+1)\cdot 2^n\cdot \frac{1-\left(\frac{-d}{2}\right)^{n+1}2^{-(n+1)} }{1-\left( \frac{-d}{2} \right)2^{-1}}$}   \\
	 	\hline
	 	{\bf odd} &  {\bf -4d} & {\bf $D^{\frac{n^2+3n}{4}}\cdot(1-2^{-(n+1)})\cdot \prod\limits_{p|d} (1+\left( \frac{(-1)^{(n+3)/2}\cdot 2}{p}\right)\cdot p^{-\frac{n+1}{2}}) \prod\limits_{j=1}^{n}\frac{j!}{(2\pi)^{j+1}}\cdot 2^n \zeta(2)  \ldots \zeta(n+1)$}   \\
	 	\hline
	 	{\bf odd} &  {\bf -d} & {\bf $D^{\frac{n^2+3n}{4}}\cdot \prod\limits_{p|d} (1+\left( \frac{(-1)^{(n+3)/2}\cdot 2}{p}\right)p^{-\frac{n+1}{2}}) \prod\limits_{j=1}^{n}\frac{j!}{(2\pi)^{j+1}}\cdot 2^n\frac{1-\left(\frac{-d}{2}\right)^{n+1}2^{-(n+1)} }{1-\left( \frac{-d}{2} \right)2^{-1}}\zeta(2)  \ldots  \zeta(n+1)$}      \\
	 	\hline
	\end{tabular}
\end{table}
\end{theorem}

\section*{Acknowledgments} The author is very grateful to O.V. Schwarzman for suggesting the problem and numerous valuable conversations. The author was supported by RSF grant no. 19-71-30002 and by Chebyshev Laboratory, St. Petersburg State University.

\section*{Bruinier formula for principal ideal rings}

The ring of integers of the field $K$ is a principal ideal ring for $$d\in  \{1, 2, 3, 7, 11, 19, 43, 67, 163\}.$$ We consider the case where $p=2$ is unramified, so $$d\in \{3, 7, 11, 19, 43, 67, 163\}.$$

\begin{definition}
Vector $e\in L$ is called a root if $(e,e)>0$ and the reflection in the subspace $e^\perp$ belongs to $\Gamma$. The length of the root $e$ is the number $(e,e)$. 
\end{definition}

The length of the roots of the lattice divides $2s$, where $s$ is the greatest invariant factor of the lattice. Since the lattice $L_n$ is unimodular, there only exist reflections in roots with length $1$ and $2$. 
\begin{definition}
A vector $l$ in a lattice $L$ is called primitive if  $K\cdot l \cap L=\mathcal{O}_K \cdot l$.
\end{definition}
Consider the number of orbits for vectors of a fixed length $a$:
$$N(L_n, a)=\{l\in L_n| (l,l)=a \textit{ and } l \textit{ is primitive}\}/Aut(L_n).$$
We use the following theorem by D. James (\cite{jam}) and formulate it in a form, convenient to us:
\begin{theorem}
Let $K=\mathbb{Q}(\sqrt{-d})$ be a quadratic field such that $p=2$ is unramified. Let $L$ be a unimodular lattice on an indefinite hermitian space with dimension $n\geq 3$. Let $a$ be a nonzero element in $\mathcal{O}_K$, represented by $L$. Then $|N(L,a)|=1$.
\end{theorem}
Since in our case $\mathcal{O}_K$ is a p.i.d., all the roots of the lattice $L_n$ are primitive. Applying James's theorem to the lattice $L_n$ we get that there is only one orbit for the roots with length $1$ or $2$.
It means, that 
\begin{itemize}
    \item the stabiliser of the root with length $1$ in the lattice $L_n$ is $U(L_{n-1}, \mathcal{O}_K)$;
    \item the stabiliser of the root with length $2$ in the lattice $L_n$ is $U(M_{n-1}, \mathcal{O}_K)$.
\end{itemize}

First we consider the case $d\neq 3$.
Each reflection has order $2$ (the group of units is $\pm 1$). 
Then the invariant $K(\Gamma)$ is
$$K(\Gamma)=\dfrac{\frac{1}{2}Vol(L_{n-1})+\frac{1}{2}Vol(M_{n-1})}{ Vol(L_n)}.$$ 
On the other hand, $K(\Gamma)\geq 2(n+1)$ (because there is $n+1$ generating function of weight at least $1$ and $n+1$ comes form the dimension).
So, if 
$$K(\Gamma)=\dfrac{\frac{1}{2}Vol(L_{n-1})+\frac{1}{2}Vol(M_{n-1})}{Vol(L_n)}\leq 2n+2,$$
then the corresponding algebra of automorphic forms can't be free.

Now we consider the case $d=3$. The difference from the previous case is that now the lattice $L_n$ is stable under the multiplication by $\epsilon E$, $\epsilon=\sqrt[6]{1}$, $E$ is the identity matrix. It follows from the definition of an automorphic form that in this case the weight of all forms is divisible by $6$. 

The reflections in roots with square $1$ have order $6$ and the reflections in roots with square $2$ have order $2$.
Then the invariant $K(\Gamma)$ is
$$K(\Gamma)=\dfrac{\frac{5}{6}Vol(L_{n-1})+\frac{1}{2}Vol(M_{n-1})}{Vol(L_n)}.$$ 
On the other hand $K\geq 7(n+1)$ (because there is $n+1$ generating function of weight at least $6$ and $n+1$ comes form the dimension).
So, if 
$$K(\Gamma)=\dfrac{\frac{5}{6}Vol(L_{n-1})+\frac{1}{2}Vol(M_{n-1})}{Vol(L_n)}\leq 7n+7,$$
then the corresponding algebra of automorphic forms can't be free.
To shorten the notation from now on we will use $K$ instead of $K(\Gamma)$.

\subsection*{The case of odd $n$}
Consider odd $n$. Then
$$K=\dfrac{\frac{1}{2}D^{\frac{(n-1)^2+3(n-1)}{4}}\cdot \prod\limits_{j=1}^{n-1}\frac{j!}{(2\pi)^{j+1}}\cdot \zeta(2)\cdot L(3)\cdot \zeta(4)\cdot L(5)\cdot \ldots \cdot L(n)\left(1+2^{n-1}\cdot \frac{1-\left(\frac{-d}{2}\right)^{n}2^{-n} }{1-\left( \frac{-d}{2} \right)2^{-1}}\right)}{D^{\frac{n^2+3n}{4}}\cdot \prod\limits_{j=1}^{n}\frac{j!}{(2\pi)^{j+1}}\cdot \zeta(2)\cdot L(3)\cdot \zeta(4)\cdot L(5)\cdot \ldots \cdot \zeta(n+1)},$$

$$K=\dfrac{(2\pi)^{n+1}\left(1+2^{n-1}\cdot \frac{1-\left(\frac{-d}{2}\right)^{n}2^{-n} }{1-\left( \frac{-d}{2} \right)2^{-1}}\right)}{2\cdot D^{(n+1)/2}\prod\limits_{p|d} (1+\left( \frac{(-1)^{(n+3)/2}}{p}\right)p^{-\frac{n+1}{2}}) \cdot n! \zeta(n+1)},$$

$$K=\dfrac{(2\pi)^{n+1}\left(1+2^{n-1}\cdot \frac{1-\left(\frac{-d}{2}\right)^{n}2^{-n} }{1-\left( \frac{-d}{2} \right)2^{-1}}\right)}{2\cdot d^{(n+1)/2}(1+\left( \frac{(-1)^{(n+3)/2}}{d}\right)d^{-\frac{n+1}{2}}) \cdot n! \zeta(n+1)},$$

$$K=\dfrac{(2\pi)^{n+1}\left(1+ \frac{2^n-\left(\frac{-d}{2}\right)^{n}}{2-\left( \frac{-d}{2} \right)}\right)}{2\left(d^{(n+1)/2}+\left( \frac{(-1)^{(n+3)/2}}{d}\right)\right) \cdot n!\cdot \zeta(n+1)}.$$

We note that $\left( \frac{-d}{2}\right)=-1$ for all $d\neq 7$ and $\left( \frac{-7}{2}\right)=1$. 

Moreover, $\left( \frac{(-1)^{(n+3)/2}}{d}\right)=1$, if $n\equiv 1 \pmod{4}$, and $\left( \frac{(-1)^{(n+3)/2}}{d}\right)=(-1)^{(d-1)/2}=-1$, if $n\equiv -1 \pmod{4}$.

We get:
\begin{table}[h]
\centering
	\begin{tabular}{|c|c|c|}
		\hline
		{\bf n} & {\bf d} & {\bf $K$} \\
	 	\hline
	 	{\bf $n\equiv 1 \pmod{4}$} &  {\bf $d\neq 7$ } & {\bf $\dfrac{(2\pi)^{n+1}\left(1+ \frac{2^n+1}{3}\right)}{2\left(d^{(n+1)/2}+1\right) \cdot n!\cdot \zeta(n+1)}$}   \\
	 	\hline
	 	{\bf $n\equiv -1 \pmod{4}$} &  {\bf $d\neq 7$} & {\bf $\dfrac{(2\pi)^{n+1}\left(1+ \frac{2^n+1}{3}\right)}{2\left(d^{(n+1)/2}-1\right) \cdot n!\cdot \zeta(n+1)}$ }   \\
	 	\hline
	 	{\bf $n\equiv 1 \pmod{4}$} &  {\bf $d=7$} & {\bf $\dfrac{(2\pi)^{n+1}\cdot 2^n}{2\left(d^{(n+1)/2}+1\right) \cdot n!\cdot \zeta(n+1)}$}   \\
	 	\hline
	 	{\bf $n\equiv -1 \pmod{4}$} &  {\bf $d=7$} & {\bf $\dfrac{(2\pi)^{n+1}\cdot 2^n}{2\left(d^{(n+1)/2}-1\right) \cdot n!\cdot \zeta(n+1)}$}      \\
	 	\hline
	\end{tabular}
\end{table}

We note that $\zeta(n+1)>1,$ so, removing this multiplier from the denominator, we get an estimate from above for $K$. Denote this by $K'$:
\newpage
\begin{table}[h]
\centering
	\begin{tabular}{|c|c|c|}
		\hline
		{\bf n} & {\bf d} & {\bf $K'$} \\
	 	\hline
	 	{\bf $n\equiv 1 \pmod{4}$} &  {\bf $d\neq 7$ } & {\bf $\dfrac{(2\pi)^{n+1}\left(1+ \frac{2^n+1}{3}\right)}{2\left(d^{(n+1)/2}+1\right) \cdot n!}$}   \\
	 	\hline
	 	{\bf $n\equiv -1 \pmod{4}$} &  {\bf $d\neq 7$} & {\bf $\dfrac{(2\pi)^{n+1}\left(1+ \frac{2^n+1}{3}\right)}{2\left(d^{(n+1)/2}-1\right) \cdot n!}$ }   \\
	 	\hline
	 	{\bf $n\equiv 1 \pmod{4}$} &  {\bf $d=7$} & {\bf $\dfrac{(2\pi)^{n+1}\cdot 2^n}{2\left(d^{(n+1)/2}+1\right) \cdot n!}$}   \\
	 	\hline
	 	{\bf $n\equiv -1 \pmod{4}$} &  {\bf $d=7$} & {\bf $\dfrac{(2\pi)^{n+1}\cdot 2^n}{2\left(d^{(n+1)/2}-1\right) \cdot n!}$}      \\
	 	\hline
	\end{tabular}
\end{table}

Since when $d$ is fixed, the numerator grows faster than the denominator, starting from some moment $K'<2n+2$.

First consider $d=7$:

When $n=3$ and $n=5$ we get $K'>2n+2$. But for $n>5$ we get $K'<2n+2$, so when $n$ is odd and $n>5$ the corresponding algebra of automorphic forms can't be free.

Let's compute the exact value of $K$ for $n=3$ and $n=5$. We use that $\zeta(4)=\frac{\pi^4}{90}$ and $\zeta(6)=\frac{\pi^6}{945}$. 
Then:
\begin{table}[h]
\centering
	\begin{tabular}{|c|c|c|}
		\hline
		{\bf n} &  {\bf $K$}& {\bf $2n+2$} \\
	 	\hline
	 	{\bf $n=5$}  & {\bf $\dfrac{(2\pi)^{6}\cdot 2^5}{2\left(7^{3}+1\right) \cdot 5!\cdot \zeta(6)}\approx 23.4$} &  {$12$}  \\
	 	\hline
	 	{\bf $n=3$} & {\bf $\dfrac{(2\pi)^{4}\cdot 2^3}{2\left(7^{2}-1\right) \cdot 3!\cdot \zeta(4)}=20$} &  {$8$}      \\
	 	\hline
	\end{tabular}
\end{table}

The value of $K$ isn't integer for $n=5$, so the corresponding algebra of automorphic forms isn't free as well.

Now we consider the case $d\geq 11$.

We note that for $d=11$ and $n=3$ or $n=5$ the inequality $K'<2n+2$ holds. When $d$ grows, the numerator of $K'$ doesn't change, and the denominator gets bigger. So, for the rest values of $d$ the inequality $K'<2n+2$ holds as well. So, the corresponding algebra of automorphic forms can't be free.

So we get that when $n$ is odd and $d\geq 7$ the only algebra of automorphic forms that can be free corresponds to the values $d=7$ and $n=3$.

\subsection*{The case of even $n$}

Writing down the Bruinier formula, we get:

$$K=\dfrac{\frac{1}{2}\zeta(2) \cdot \ldots \cdot L(n)\left(1+\left( \frac{(-1)^{(n+2)/2}}{d}\right)d^{-\frac{n}{2}} + 
(1+\left( \frac{(-1)^{(n+2)/2}\cdot 2}{d}\right)d^{-\frac{n}{2}})  2^{n-1}\frac{1-\left(\frac{-d}{2}\right)^{n}\cdot 2^{-n} }{1-\left( \frac{-d}{2} \right)2^{-1}}
\right)}{d^{\frac{n+1}{2}} \cdot\frac{n!}{(2\pi)^{n+1}}\cdot \zeta(2)\cdot L(3)\cdot  \ldots \cdot L(n+1)},$$

$$K=\dfrac{(2\pi)^{n+1}\left(1+\left( \frac{(-1)^{(n+2)/2}}{d}\right)d^{-\frac{n}{2}} + 
(1+\left( \frac{(-1)^{(n+2)/2}\cdot 2}{d}\right)d^{-\frac{n}{2}})  2^{n-1}\frac{1-\left(\frac{-d}{2}\right)^{n}\cdot 2^{-n} }{1-\left( \frac{-d}{2} \right)2^{-1}}
\right)}{2d^{\frac{n+1}{2}} \cdot n!\cdot L(n+1)},$$

$$K=\dfrac{(2\pi)^{n+1}\left(1+\left( \frac{(-1)^{(n+2)/2}}{d}\right)d^{-\frac{n}{2}} + 
(1+\left( \frac{(-1)^{(n+2)/2}\cdot 2}{d}\right)d^{-\frac{n}{2}})  \frac{2^n-\left(\frac{-d}{2}\right)^{n} }{2-\left( \frac{-d}{2} \right)}
\right)}{2d^{\frac{n+1}{2}} \cdot n!\cdot L(n+1)}.$$

We note that $\left( \frac{-d}{2}\right)=-1$ for all $d\neq 7$ and $\left( \frac{-7}{2}\right)=1$. So:

$$K=\dfrac{(2\pi)^{n+1}\left(1+\left( \frac{(-1)^{(n+2)/2}}{d}\right)d^{-\frac{n}{2}} + 
(1+\left( \frac{(-1)^{(n+2)/2}\cdot 2}{d}\right)d^{-\frac{n}{2}})  \frac{2^n-1 }{3}
\right)}{2d^{\frac{n+1}{2}} \cdot n!\cdot L(n+1)}, ~d\neq 7,$$
$$K=\dfrac{(2\pi)^{n+1}\left(1+\left( \frac{(-1)^{(n+2)/2}}{d}\right)d^{-\frac{n}{2}} + 
(1+\left( \frac{(-1)^{(n+2)/2}\cdot 2}{d}\right)d^{-\frac{n}{2}})  (2^n-1) 
\right)}{2d^{\frac{n+1}{2}} \cdot n!\cdot L(n+1)}, ~d=7.$$

Since $\left( \frac{-1}{d} \right)=(-1)^{\frac{d-1}{2}}=-1$, $\left(\frac{2}{d}\right)=(-1)^{\frac{d^2-1}{8}}=-1$ for $d\neq 7$ and $\left( \frac{2}{7} \right)=1$ and $\left( \frac{-2}{d}\right)=\left(\frac{-1}{d}\right)\cdot \left( \frac{2}{d}\right)$, we get the following table:

\begin{table}[h]
\centering
	\begin{tabular}{|c|c|c|}
		\hline
		{\bf n} & {\bf d} & {\bf $K$} \\
	 	\hline
	 	{\bf $n\equiv 0 \pmod{4}$} &  {\bf $d\neq 7$ } & {\bf $\dfrac{(2\pi)^{n+1}\left(1- d^{-n/2}+\frac{2^n-1}{3}(1+d^{-n/2})\right)}{2\cdot d^{(n+1)/2 }\cdot n!\cdot L(n+1)}$}   \\
	 	\hline
	 	{\bf $n\equiv 2 \pmod{4}$} &  {\bf $d\neq 7$} & {\bf $\dfrac{(2\pi)^{n+1}\left(1+ d^{-n/2}+\frac{2^n-1}{3}(1-d^{-n/2})\right)}{2\cdot d^{(n+1)/2 }\cdot n!\cdot L(n+1)}$ }   \\
	 	\hline
	 	{\bf $n\equiv 0 \pmod{4}$} &  {\bf $d=7$} & {\bf $\dfrac{(2\pi)^{n+1}\left(1- d^{-n/2}+(2^n-1)(1-d^{-n/2})\right)}{2\cdot d^{(n+1)/2 }\cdot n!\cdot L(n+1)}$}   \\
	 	\hline
	 	{\bf $n\equiv 2 \pmod{4}$} &  {\bf $d=7$} & {\bf $\dfrac{(2\pi)^{n+1}\left(1+ d^{-n/2}+(2^n-1)(1+d^{-n/2})\right)}{2\cdot d^{(n+1)/2 }\cdot n!\cdot L(n+1)}$}      \\
	 	\hline
	\end{tabular}
\end{table}

Simplifying the value of $K$ for $d=7$, we get:
\begin{table}[h]
\centering
	\begin{tabular}{|c|c|c|}
		\hline
		{\bf n} & {\bf d} & {\bf $K$} \\
	 	\hline
	 	{\bf $n\equiv 0 \pmod{4}$} &  {\bf $d\neq 7$ } & {\bf $\dfrac{(2\pi)^{n+1}\left(1- d^{-n/2}+\frac{2^n-1}{3}(1+d^{-n/2})\right)}{2\cdot d^{(n+1)/2 }\cdot n!\cdot L(n+1)}$}   \\
	 	\hline
	 	{\bf $n\equiv 2 \pmod{4}$} &  {\bf $d\neq 7$} & {\bf $\dfrac{(2\pi)^{n+1}\left(1+ d^{-n/2}+\frac{2^n-1}{3}(1-d^{-n/2})\right)}{2\cdot d^{(n+1)/2 }\cdot n!\cdot L(n+1)}$ }   \\
	 	\hline
	 	{\bf $n\equiv 0 \pmod{4}$} &  {\bf $d=7$} & {\bf $\dfrac{\pi^{n+1}\cdot2^{2n}(1-d^{-n/2})}{d^{(n+1)/2 }\cdot n!\cdot L(n+1)}$}   \\
	 	\hline
	 	{\bf $n\equiv 2 \pmod{4}$} &  {\bf $d=7$} & {\bf $\dfrac{\pi^{n+1}\cdot2^{2n}(1+d^{-n/2})}{d^{(n+1)/2 }\cdot n!\cdot L(n+1)}$}      \\
	 	\hline
	\end{tabular}
\end{table}

We note that $\frac{1}{L(n+1)}\leq 2$. This gives the estimate from the above for $K$. Denote this by $K'$:
\newpage
\begin{table}[h]
\centering
	\begin{tabular}{|c|c|c|}
		\hline
		{\bf n} & {\bf d} & {\bf $K'$} \\
	 	\hline
	 	{\bf $n\equiv 0 \pmod{4}$} &  {\bf $d\neq 7$ } & {\bf $\dfrac{(2\pi)^{n+1}\left(1- d^{-n/2}+\frac{2^n-1}{3}(1+d^{-n/2})\right)}{d^{(n+1)/2 }\cdot n!}$}   \\
	 	\hline
	 	{\bf $n\equiv 2 \pmod{4}$} &  {\bf $d\neq 7$} & {\bf $\dfrac{(2\pi)^{n+1}\left(1+ d^{-n/2}+\frac{2^n-1}{3}(1-d^{-n/2})\right)}{d^{(n+1)/2 }\cdot n!}$ }   \\
	 	\hline
	 	{\bf $n\equiv 0 \pmod{4}$} &  {\bf $d=7$} & {\bf $\dfrac{\pi^{n+1}\cdot2^{2n+1}(1-d^{-n/2})}{d^{(n+1)/2 }\cdot n!}$}   \\
	 	\hline
	 	{\bf $n\equiv 2 \pmod{4}$} &  {\bf $d=7$} & {\bf $\dfrac{\pi^{n+1}\cdot2^{2n
	 	+1}(1+d^{-n/2})}{d^{(n+1)/2 }\cdot n!}$}      \\
	 	\hline
	\end{tabular}
\end{table} 

We estimate the numerator for $d\neq 7$ from the above, using that $d^{-n/2}<\frac{1}{2},~-d^{n/2}<0$. Denote this number by $K''$.

\begin{table}[h]
\centering
	\begin{tabular}{|c|c|c|}
		\hline
		{\bf n} & {\bf d} & {\bf $K''$} \\
	 	\hline
	 	{\bf $n\equiv 0 \pmod{4}$} &  {\bf $d\neq 7$ } & {\bf $\dfrac{(2\pi)^{n+1}\left(1+\frac{2^n-1}{3}(1+\frac{1}{2})\right)}{d^{(n+1)/2 }\cdot n!}$}   \\
	 	\hline
	 	{\bf $n\equiv 2 \pmod{4}$} &  {\bf $d\neq 7$} & {\bf $\dfrac{(2\pi)^{n+1}\left(1+ \frac{1}{2}+\frac{2^n-1}{3}\cdot 1\right)}{d^{(n+1)/2 }\cdot n!}$ }   \\
	 	\hline
	 	{\bf $n\equiv 0 \pmod{4}$} &  {\bf $d=7$} & {\bf $\dfrac{\pi^{n+1}\cdot2^{2n+1}(1-d^{-n/2})}{d^{(n+1)/2 }\cdot n!}$}   \\
	 	\hline
	 	{\bf $n\equiv 2 \pmod{4}$} &  {\bf $d=7$} & {\bf $\dfrac{\pi^{n+1}\cdot2^{2n
	 	+1}(1+d^{-n/2})}{d^{(n+1)/2 }\cdot n!}$}      \\
	 	\hline
	\end{tabular}
\end{table} 

Since when $d$ is fixed, the numerator grows faster than the denominator, starting from some moment $K''<2n+2$.

First consider the case $d=7$:

When $n=2$, $n=4$ and $n=6$ we get the inequality $K''>2n+2$. But for even $n>6$ the inequality $K''<2n+2$ holds, so the corresponding algebra of automorphic forms can't be free.

We compute the exact values of $K$ for $n=2$, $n=4$ and $n=6$, using that for the field $\mathbb{Q}(\sqrt{-7})$ the values of $L$-function are: $L(3)=\frac{32}{2401}\sqrt{7}\pi^3$, $L(5)=\frac{64}{50421}\sqrt{7}\pi^5$, $L(7)=\frac{4672}{37059435}\sqrt{7}\pi^7$. Then:
\begin{table}[h]
\centering
	\begin{tabular}{|c|c|c|}
		\hline
		{\bf n} &  {\bf $K$}& {\bf $2n+2$} \\
	 	\hline
	 	{\bf $n=2$}  & {\bf $\dfrac{\pi^{3}\cdot2^{4}(1+7^{-1})}{7^{3/2 }\cdot 2!\cdot L(3)}=14$} &  {$6$}  \\
	 	\hline
	 	{\bf $n=4$} & {\bf $\dfrac{\pi^{5}\cdot2^{8}(1-7^{-2})}{7^{5/2 }\cdot 4!\cdot L(5)}=24$} &  {$10$}      \\
	 	\hline
	 	{\bf $n=6$} & {\bf $\dfrac{\pi^{7}\cdot2^{12}(1+7^{-3})}{7^{7/2 }\cdot 6!\cdot L(7)}\approx 18.8$} &  {$12$}      \\
	 	\hline
	\end{tabular}
\end{table}

The value of $K$ is not integer for $n=6$, so the corresponding algebra can't be free.

Now we consider the case $d\geq 11$.

We note that for $d=11$ and $n=4$ or $n=6$ the inequality $K''<2n+2$ holds. When $d$ grows, the numerator of $K''$ doesn't change, and the denominator gets bigger. So, for the rest values of $d$ the inequality $K''<2n+2$ holds as well. So, the corresponding algebra of automorphic forms can't be free.

%So we get that when $n$ is odd and $d\geq 7$ the only algebra of automorphic forms that can be free corresponds to the values $d=7$ and $n=3$. 

We compute the exact value of $K$ for $n=2$ and $d=11$, $d=19$, $d=43$, $d=67$, $d=163$. We use the exact values of $L$-functions: $L(3)_{11}=\frac{12}{1331}\sqrt{11}\pi^3$, $L(3)_{19}=\frac{44}{6859}\sqrt{19}\pi^3$, $L(3)_{43}=\frac{332}{79507}\sqrt{43}\pi^3$, $L(3)_{67}=\frac{1004}{300763}\sqrt{67}\pi^3$, $L(3)_{163}=\frac{9260}{4330747}\sqrt{163}\pi^3$($L(3)_d$ is $L(3)$ for the field $\mathbb{Q}(\sqrt{-d})$). 

Since $n=2$ then $2n+2=6$.
We get:
\begin{table}[h]
\centering
	\begin{tabular}{|c|c|}
		\hline
		{\bf d} &  {\bf $K$} \\
	 	\hline
	 	{\bf $d=11$}  & {\bf $\dfrac{(2\pi)^{3}}{11^{3/2 }\cdot 2\cdot L(3)}=\frac{11}{3}<6$}   \\
	 	\hline
	 	{\bf $d=19$} & {\bf $\dfrac{(2\pi)^{3}}{19^{3/2 }\cdot 2\cdot L(3)}=\frac{19}{11}<6$}       \\
	 	\hline
	 	{\bf $d=43$} & {\bf $\dfrac{(2\pi)^{3}}{43^{3/2 }\cdot 2\cdot L(3)=\frac{43}{83}<6}$}     \\
	 	\hline
	 	{\bf $d=67$} & {\bf $\dfrac{(2\pi)^{3}}{67^{3/2 }\cdot 2\cdot L(3)=\frac{67}{251}<6}$}       \\
	 	\hline
	 	{\bf $d=163$} & {\bf $\dfrac{(2\pi)^{3}}{163^{3/2 }\cdot 2\cdot L(3)}=\frac{163}{2315}<6$}      \\
	 	\hline
	\end{tabular}
\end{table}

It follows that in all this cases the corresponding algebra of automorphic forms can't be free.

We get that for even $n$ and $d\geq 7$ the algebra of automorphic forms which can be free corresponds to the values $d=7$ and $n=2$ or $d=7$ and $n=4$.

\section*{$d=3$}
Here we provide the same calculations as in the previous case, just the Bruinier formula looks a bit different:
$$\frac{5}{6}Vol(L_{n-1})+\frac{1}{2}Vol(M_{n-1})=K\cdot Vol(L_n),$$ where $K\geq 7(n+1)$ (because there is $n+1$ generating function of weight at least $6$ and $n+1$ comes form the dimension).
So, if 
$$K=\dfrac{\frac{5}{6}Vol(L_{n-1})+\frac{1}{2}Vol(M_{n-1})}{Vol(L_n)}\leq 7n+7,$$
then the corresponding algebra of automorphic forms can't be free.

\subsection*{The case of odd $n$}

Consider odd $n$. Then
$$K=\dfrac{\frac{1}{2}D^{\frac{(n-1)^2+3(n-1)}{4}}\cdot \prod\limits_{j=1}^{n-1}\frac{j!}{(2\pi)^{j+1}}\cdot \zeta(2)\cdot L(3)\cdot \zeta(4)\cdot L(5)\cdot \ldots \cdot L(n)\left(\frac{5}{3}+2^{n-1}\cdot \frac{1-\left(\frac{-d}{2}\right)^{n}2^{-n} }{1-\left( \frac{-d}{2} \right)2^{-1}}\right)}{D^{\frac{n^2+3n}{4}}\cdot \prod\limits_{j=1}^{n}\frac{j!}{(2\pi)^{j+1}}\cdot \zeta(2)\cdot L(3)\cdot \zeta(4)\cdot L(5)\cdot \ldots \cdot \zeta(n+1)},$$

$$K=\dfrac{(2\pi)^{n+1}\left(\frac{5}{3}+2^{n-1}\cdot \frac{1-\left(\frac{-3}{2}\right)^{n}2^{-n} }{1-\left( \frac{-3}{2} \right)2^{-1}}\right)}{2\cdot 3^{(n+1)/2}(1+\left( \frac{(-1)^{(n+3)/2}}{3}\right)3^{-\frac{n+1}{2}}) \cdot n! \zeta(n+1)},$$

$$K=\dfrac{(2\pi)^{n+1}\left(\frac{5}{3}+2^{n-1}\cdot \frac{1-\left(\frac{-3}{2}\right)^{n}2^{-n} }{1-\left( \frac{-3}{2} \right)2^{-1}}\right)}{2\cdot 3^{(n+1)/2}(1+\left( \frac{(-1)^{(n+3)/2}}{3}\right)3^{-\frac{n+1}{2}}) \cdot n! \zeta(n+1)},$$

$$K=\dfrac{(2\pi)^{n+1}\left(\frac{5}{3}+ \frac{2^n-\left(\frac{-3}{2}\right)^{n}}{2-\left( \frac{-3}{2} \right)}\right)}{2\left(3^{(n+1)/2}+\left( \frac{(-1)^{(n+3)/2}}{3}\right)\right) \cdot n!\cdot \zeta(n+1)}.$$

We note that $\left( \frac{-3}{2}\right)=-1$. Moreover, $\left( \frac{(-1)^{(n+3)/2}}{3}\right)=1$, if $n\equiv 1 \pmod{4}$, and $\left( \frac{(-1)^{(n+3)/2}}{3}\right)=-1$, if $n\equiv -1 \pmod{4}$.

We get:
\begin{table}[h]
\centering
	\begin{tabular}{|c|c|}
		\hline
		{\bf n}  & {\bf $K$} \\
	 	\hline
	 	{\bf $n\equiv 1 \pmod{4}$} & {\bf $\dfrac{(2\pi)^{n+1}\left(6+ 2^n\right)}{6\left(3^{(n+1)/2}+1\right) \cdot n!\cdot \zeta(n+1)}$}   \\
	 	\hline
	 	{\bf $n\equiv -1 \pmod{4}$}  & {\bf $\dfrac{(2\pi)^{n+1}\left(6+ 2^n\right)}{6\left(3^{(n+1)/2}-1\right) \cdot n!\cdot \zeta(n+1)}$ }   \\
	 	\hline
	\end{tabular}
\end{table}

We note that $\zeta(n+1)>1,$ so, removing this multiplier from the denominator, we get an estimate from above for $K$. Denote this by $K'$:
\begin{table}[h]
\centering
	\begin{tabular}{|c|c|}
		\hline
		{\bf n}  & {\bf $K$} \\
	 	\hline
	 	{\bf $n\equiv 1 \pmod{4}$} & {\bf $\dfrac{(2\pi)^{n+1}\left(6+ 2^n\right)}{6\left(3^{(n+1)/2}+1\right) \cdot n!}$}   \\
	 	\hline
	 	{\bf $n\equiv -1 \pmod{4}$}  & {\bf $\dfrac{(2\pi)^{n+1}\left(6+ 2^n\right)}{6\left(3^{(n+1)/2}-1\right) \cdot n!}$ }   \\
	 	\hline
	\end{tabular}
\end{table}

Since the numerator grows faster than the denominator, starting from some moment $K'<7n+7$.

For odd $n>9$ we get $K'<7n+7$, so when $n$ is odd and $n>9$ the corresponding algebra of automorphic forms can't be free.

Let's compute the exact value of $K$ for $n=3$, $n=5$, $n=7$ and $n=9$. We use that $\zeta(4)=\frac{\pi^4}{90}$, $\zeta(6)=\frac{\pi^6}{945}$, $\zeta(8)=\frac{\pi^8}{9450}$ and $\zeta(10)=\frac{\pi^{10}}{93555}$. 
Then:
\begin{table}[h]
\centering
	\begin{tabular}{|c|c|c|}
		\hline
		{\bf n} &  {\bf $K$}& {\bf $7n+7$} \\
	 	\hline
	 	{\bf $n=9$}  & {\bf $\frac{5698}{61}$} &  {$70$}  \\
	 	\hline
	 	{\bf $n=7$}  & {\bf $134$} &  {$56$}  \\
	 	\hline
	 	{\bf $n=5$}  & {\bf $114$} &  {$42$}  \\
	 	\hline
	 	{\bf $n=3$} & {\bf 70} &  {$28$}      \\
	 	\hline
	\end{tabular}
\end{table}

The value of $K$ isn't integer for $n=9$, so the corresponding algebra of automorphic forms isn't free as well. 

So we get that when $n$ is odd and $d=3$ the algebras of automorphic forms that can be free correspond to the values $n=7$, $n=5$ and $n=3$.

\subsection*{The case of even $n$}

Writing down the Bruinier formula, we get:

$$K=\dfrac{\frac{1}{2}\zeta(2) \cdot \ldots \cdot \zeta(n)\left(\frac{5}{3}(1+\left( \frac{(-1)^{(n+2)/2}}{d}\right)d^{-\frac{n}{2}}) + 
(1+\left( \frac{(-1)^{(n+2)/2}\cdot 2}{d}\right)d^{-\frac{n}{2}})  2^{n-1}\frac{1-\left(\frac{-d}{2}\right)^{n}\cdot 2^{-n} }{1-\left( \frac{-d}{2} \right)2^{-1}}
\right)}{d^{\frac{n+1}{2}} \cdot\frac{n!}{(2\pi)^{n+1}}\cdot \zeta(2)\cdot L(3)\cdot  \ldots \cdot L(n+1)},$$

$$K=\dfrac{(2\pi)^{n+1}\left(\frac{5}{3}(1+\left( \frac{(-1)^{(n+2)/2}}{3}\right)3^{-\frac{n}{2}}) + 
(1+\left( \frac{(-1)^{(n+2)/2}\cdot 2}{3}\right)3^{-\frac{n}{2}})  2^{n-1}\frac{1-\left(\frac{-3}{2}\right)^{n}\cdot 2^{-n} }{1-\left( \frac{-3}{2} \right)2^{-1}}
\right)}{2\cdot3^{\frac{n+1}{2}} \cdot n!\cdot L(n+1)},$$

$$K=\dfrac{(2\pi)^{n+1}\left(\frac{5}{3}(1+\left( \frac{(-1)^{(n+2)/2}}{3}\right)3^{-\frac{n}{2}}) + 
(1+\left( \frac{(-1)^{(n+2)/2}\cdot 2}{3}\right)3^{-\frac{n}{2}})  \frac{2^n-\left(\frac{-3}{2}\right)^{n} }{2-\left( \frac{-3}{2} \right)}
\right)}{2\cdot 3^{\frac{n+1}{2}} \cdot n!\cdot L(n+1)}.$$

We note that $\left( \frac{-3}{2}\right)=-1$. So:

$$K=\dfrac{(2\pi)^{n+1}\left(\frac{5}{3}(1+\left( \frac{(-1)^{(n+2)/2}}{3}\right)3^{-\frac{n}{2}} )+ 
(1+\left( \frac{(-1)^{(n+2)/2}\cdot 2}{3}\right)3^{-\frac{n}{2}})  \frac{2^n-1 }{3}
\right)}{2\cdot 3^{\frac{n+1}{2}} \cdot n!\cdot L(n+1)}.$$

Since $\left( \frac{-1}{3} \right)=\left(\frac{2}{3}\right)=-1$ and $\left( \frac{-2}{3}\right)=1$, we get the following table:

\begin{table}[h]
\centering
	\begin{tabular}{|c|c|}
		\hline
		{\bf n} & {\bf $K$} \\
	 	\hline
	 	{\bf $n\equiv 0 \pmod{4}$} & {\bf $ \dfrac{(2\pi)^{n+1}\left(\frac{5}{3}(1-3^{-\frac{n}{2}} )+ 
(1+3^{-\frac{n}{2}})  \frac{2^n-1 }{3}
\right)}{2\cdot 3^{\frac{n+1}{2}} \cdot n!\cdot L(n+1)}.$}   \\
	 	\hline
	 	{\bf $n\equiv 2 \pmod{4}$} & {\bf $\dfrac{(2\pi)^{n+1}\left(\frac{5}{3}(1+3^{-\frac{n}{2}} )+ 
(1-3^{-\frac{n}{2}})  \frac{2^n-1 }{3}
\right)}{2\cdot 3^{\frac{n+1}{2}} \cdot n!\cdot L(n+1)}.$ }   \\
	 	\hline
	\end{tabular}
\end{table}

Simplifying the value of $K$, we get:
\begin{table}[h]
\centering
	\begin{tabular}{|c|c|}
		\hline
		{\bf n} & {\bf $K$} \\
	 	\hline
	 	{\bf $n\equiv 0 \pmod{4}$} & {\bf $ \dfrac{(2\pi)^{n+1}\left(5(1-3^{-\frac{n}{2}} )+ 
(1+3^{-\frac{n}{2}})(2^n-1)
\right)}{6\cdot 3^{\frac{n+1}{2}} \cdot n!\cdot L(n+1)}.$}   \\
	 	\hline
	 	{\bf $n\equiv 2 \pmod{4}$} & {\bf $\dfrac{(2\pi)^{n+1}\left(5(1+3^{-\frac{n}{2}} )+ 
(1-3^{-\frac{n}{2}})(2^n-1)
\right)}{6\cdot 3^{\frac{n+1}{2}} \cdot n!\cdot L(n+1)}.$ }   \\
	 	\hline
	\end{tabular}
\end{table}

We note that $\frac{1}{L(n+1)}\leq 2$. This gives the estimate from the above for $K$. Denote this by $K'$:
\newpage
\begin{table}[h]
\centering
	\begin{tabular}{|c|c|}
		\hline
		{\bf n} & {\bf $K'$} \\
	 	\hline
	 	{\bf $n\equiv 0 \pmod{4}$} & {\bf $ \dfrac{(2\pi)^{n+1}\left(5(1-3^{-\frac{n}{2}} )+ 
(1+3^{-\frac{n}{2}})(2^n-1)
\right)}{3\cdot 3^{\frac{n+1}{2}} \cdot n!}.$}   \\
	 	\hline
	 	{\bf $n\equiv 2 \pmod{4}$} & {\bf $\dfrac{(2\pi)^{n+1}\left(5(1+3^{-\frac{n}{2}} )+ 
(1-3^{-\frac{n}{2}})(2^n-1)
\right)}{3\cdot 3^{\frac{n+1}{2}} \cdot n!}.$ }   \\
	 	\hline
	\end{tabular}
\end{table}

Since the numerator grows faster than the denominator, starting from some moment $K'<7n+7$.
For odd $n>10$ we get $K'<7n+7$, so when $n$ is even and $n>10$ the corresponding algebra of automorphic forms can't be free.

Let's compute the exact value of $K$ for $n=2$, $n=4$, $n=6$, $n=8$ and $n=10$. We use that for the field $\mathbb{Q}(\sqrt{-3})$ the values of $L$-function are: $L(3)=\frac{4}{243}\sqrt{3}\pi^3$, $L(5)=\frac{4}{2187}\sqrt{3}\pi^5$, $L(7)=\frac{56}{295245}\sqrt{3}\pi^7$, $L(9)=\frac{3236}{167403915}\sqrt{3}\pi^{9}$, $L(11)=\frac{14776}{7533176175}\sqrt{3}\pi^{11}$. Then:
\begin{table}[h]
\centering
	\begin{tabular}{|c|c|c|}
		\hline
		{\bf n} &  {\bf $K$}& {\bf $7n+7$} \\
	 	\hline
	 	{\bf $n=2$}  & {\bf $39$} &  {$21$}  \\
	 	\hline
	 	{\bf $n=4$} & {\bf $95$} &  {$35$}      \\
	 	\hline
	 	{\bf $n=6$} & {\bf $127$} &  {$49$}      \\
	 	\hline
	 	{\bf $n=8$} & {\bf $\approx 118.5$} &  {$63$}      \\
	 	\hline
	 	{\bf $n=10$} & {\bf $\approx 67.3$} &  {$77$}      \\
	 	\hline
	\end{tabular}
\end{table}

The value of $K$ is not integer for $n=8$ and $n=10$, so the corresponding algebra can't be free. So we get that when $n$ is even and $d=3$ the algebras of automorphic forms that can be free correspond to the values $n\leq 7$.

\end{document}